\documentstyle[12pt]{article}
\newtheorem{theorem}{Theorem}
\newtheorem{lemma}[theorem]{Lemma}

\newtheorem{problem}[theorem]{Problem}

\newcommand{\text}[1]{\quad\mbox{#1}\quad}
\def\beq{\begin{equation}}\def\eeq{\end{equation}}
\def\beqn{\begin{eqnarray}}\def\eeqn{\end{eqnarray}}
\def\pont{\hspace{-6pt}{\bf.\ }}

\def\qed{\ifhmode\unskip\nobreak\fi\quad\ifmmode\Box\else$\Box$\fi}

\title{Coloring $2$-intersecting hypergraphs}

\author{
Lucas Colucci\\
\small Instituto de Matemática e Estatística\\[-0.8ex]
\small Universidade de São Paulo\\[-0.8ex]
\small Rua do Matão 1010, 05508-090\\[-0.8ex]
\small São Paulo, Brazil\\ \small
\texttt{lucas.colucci.souza@gmail.com}\and
Andr\'as Gy\'arf\'as\thanks{Research coordinator of the junior author at the Elective Undergraduate Research Program of Budapest Semesters in Mathematics, 2013 Summer program}\\
\small Computer and Automation Research Institute\\[-0.8ex]
\small Hungarian Academy of Sciences\\[-0.8ex]
\small Budapest, P.O. Box 63\\[-0.8ex]
\small Budapest, Hungary, H-1518\\\small
\texttt{gyarfas.andras@renyi.mta.hu}}

\begin{document}
\maketitle

\begin{abstract} A hypergraph is $2$-intersecting if any two edges intersect in at least two vertices.
Blais, Weinstein and Yoshida asked (as a first step to a more general problem) whether every $2$-intersecting hypergraph has a vertex coloring with a constant number of colors so that each hyperedge has at least $\min\{|e|,3\}$ colors. We show that there is such a coloring with at most 5 colors (which is best possible).  \end{abstract}

A {\em proper coloring} of a hypergraph is a coloring of its vertices so that no edge is monochromatic, i.e. contains at least two vertices with distinct colors.  It is  well-known that intersecting hypergraphs without singleton edges have proper colorings with at most three colors. This statement is from the seminal paper of Erd\H os and Lov\'asz \cite{EL}. Recently Blais, Weinstein and Yoshida suggested a generalization in \cite{BWY}. They consider $t$-intersecting hypergraphs, in which any two edges intersect in at least $t$ vertices and they call a coloring of the vertices {\em $c$-strong} if every edge $e$ is colored with at least $\min\{|e|,c\}$ distinct colors. One of the problems they consider is the following.

\begin{problem}\pont \label{BWYpr} (\cite{BWY}) Suppose that $\cal{H}$ is a  $t$-intersecting hypergraph. Is there a $(t+1)$-strong vertex coloring of $\cal{H}$ where the number of colors is bounded by a function of $t$? In particular, is there a $t+1$-strong vertex coloring with at most $2t+1$ colors? If true, it would be best possible, as the $2t$-element sets of a $3t$ element set demonstrate.
\end{problem}

Notice that for $t=1$ the answer to Problem \ref{BWYpr} is affirmative (for both parts) according to the starting remark but open for $t\ge 2$  \cite{BWY}. Our aim is to give an affirmative answer to both parts of the problem in case of $t=2$. Notice that intersecting hypergraphs do not always have $3$-strong colorings with any fixed number of colors: if every edge of a $(k+1)$-chromatic graph is extended by the same new vertex, the resulting intersecting hypergraph has no $3$-strong coloring with $k$ colors. Thus the $2$-intersecting condition is important in the following theorem.

\begin{theorem}\pont\label{2intersect} Every $2$-intersecting hypergraph $G$ has a $3$-strong coloring with at most five colors.
\end{theorem}

We also prove a lemma that will be used in the proof of Theorem \ref{2intersect} but has independent interest.
A hypergraph has property $P_t$ for some integer $t\ge 2$ if any $i$ edges intersect in at least $t+1-i$ vertices, for all $i, 2\le i \le t$.

\begin{lemma}\pont\label{str}Suppose that $\cal{H}$ is a hypergraph with property $P_t$. Then $\cal{H}$ has a $t$-strong coloring with at most $t+1$ colors.
\end{lemma}

\noindent {\bf Proof. \rm} Let $\cal{H}$ be a hypergraph with property $P_t$ for $t\ge 2$. Select an edge $e_0$ of $\cal{H}$ which is minimal for containment. Let $\cal{F}$ be the hypergraph defined on the vertex set of $e_0$ with edge set $\{h\cap e_0: h\in E({\cal{H}})\}$. Color each vertex not in $e$ with color $t+1$. If $t=2$, color the vertices of $e$ arbitrarily using colors 1,2 (or just color 1 if $e$ has just one vertex). Otherwise, since $\cal{F}$ has property $P_{t-1}$, we can find by induction a $(t-1)$-strong coloring $C$ on $\cal{F}$ with colors $1,2,\dots, t$. Since for each edge $h\in {\cal{H}}, |h\cap e_0|\ge t-1$, $C$ uses at least $t-1$ colors on $h\cap e_0$ and $h$ also has at least one vertex of color $t+1$. Therefore we have a $t$-strong coloring of $\cal{H}$ with $t+1$ colors. \qed

It is worth noting that Lemma \ref{str} does not hold if we require a $t$-strong coloring with at most $t$ colors. Indeed, all $t$-sets of $t+1$ elements have property $P_t$ but a $t$-strong coloring must use $t+1$ colors.
\bigskip

\noindent {\bf Proof of Theorem \ref{2intersect}.} By the condition, there are no singleton edges.
Also, a $3$-strong coloring on the minimal edges of $G$ is also a $3$-strong
coloring on $G$, thus we may assume that $G$ is an antichain.

If any three edges of $G$ have non-empty
intersection, we can apply Lemma \ref{str} and get a $3$-strong coloring with at most $4$ colors.
Thus, we may suppose that $G$ contains three edges with empty intersection, select them
with the smallest possible union, let these edges be $e_1,e_2,e_3$ and set $X=e_1\cup e_2\cup e_3$.
A vertex $v\in X$ is called a private part of $e_i$ ($i=1,2,3$) if $v\in e_i$ but $v$ is not covered by any of the other two $e_j$-s.

We color the vertices in $X$ as follows. The private parts of $e_1,e_2,e_3$ (if they exist) are colored with $1,2,3$
respectively . Notice that each intersection has at
least two vertices, color $e_1\cap e_3$ with colors $1,3$ so that
color $1$ is used only once, color $e_1\cap e_2$ with colors $2,4$
so that color $2$ is used only once. Vertices in $e_2\cap e_3$ are
all colored with color $5$.

The coloring outside $X$ varies according to the number of private parts of $e_i$-s.

\noindent {\bf Case 1. } Each $e_i$ has private parts, $i=1,2,3$.

\noindent Here we color vertices not covered by $X$ one-by one with $1$ or $2$ by the
following greedy type algorithm: if an uncolored vertex $w\notin X$
completes an edge $f$ such that all vertices of $f-\{w\}$ are colored with colors
$2,3$ only (not necessarily with both) then color $w$ with color
$1$, otherwise color it with color $2$. We claim that a $3$-strong
coloring is obtained.

Suppose there is an edge $f_{ij}$ with colors $i,j$ only, $1\le i <
j\le 5$. Edges $f_{12},f_{14},f_{24}$ would intersect $e_3$ in at
most one vertex, edge $f_{25}$ would intersect $e_1$ in at most one
vertex  and $f_{13}$ would not intersect $e_2$ at all.
Edges $f_{35},f_{45}$ would form a proper subset of $e_3,e_2$, respectively,
contradicting the antichain property.

Edge $f_{34}$ cannot exist because  the triple $f_{34},e_2,e_3$
 has no intersection and $Y=f_{34}\cup  e_2\cup e_3$  is a proper
subset of $X$ because $e_1$ has a private vertex.  Thus we get
a contradiction with the definition of $e_1,e_2,e_3$. The same
argument can be applied to exclude $f_{15},f_{23}\subset X$ (with $Y=f_{15}\cup  e_1\cup e_2, Y=f_{23}\cup  e_2\cup e_3$ and using that $e_3,e_1$ have  private vertices).

Thus the only possibility is that there is an edge $f_{15}$ or
$f_{23}$ with some vertex $w\notin X$. However, no such $f_{15}$
exists since $w\notin X$ is colored with $1$ only if there exists
edge $f$ of $G$ such that $f-\{w\}$ is colored with colors $2,3$
only thus $|f\cap f_{15}|=1$ contradiction. Moreover, no such
$f_{23}$ can exist either, because its vertex in $V-X$ colored last got color $1$ according to the rule governing Case 1.

\noindent {\bf Case 2.} Two of $e_1,e_2,e_3$ have private parts, by suitable relabeling we may suppose that the private part of $e_2$ is empty.

\noindent In this case vertices not covered by $X$ are colored with color $2$ and claim that we have a $3$-strong coloring.
The nonexistence of $f_{12},f_{13},f_{14},f_{24},f_{25}$ follow as in Case 1 and here $f_{23}$ can be excluded the same way since $|f_{23}\cap e_2|\le 1$. The exclusion of $f_{34},f_{35},f_{45}$ and $f_{15}\subset X$ is also exactly the same as in Case 1. Thus here we have to exclude only  the existence of an edge $f_{15}$ containing some vertices $w\notin X$. However, this cannot happen since here every vertex outside $X$ is colored with color $2$.

\noindent {\bf Case 3.} Exactly one of $e_1,e_2,e_3$ has a private part, by suitable relabeling we may suppose that it is $e_2$.

\noindent Here all vertices not covered by $X$ are colored with $1$.  Edges $f_{12},f_{13},f_{14},f_{15},f_{24},f_{25}$ are all excluded since there is some $e_i$ intersecting them in at most one vertex. The edges $f_{34},f_{35},f_{45}$ are excluded since they are proper subsets of some $e_i$.
The only possible edge is $f_{23}$ but in this case we can replace the triple $e_1,e_2,e_3$ by the non-intersecting triple $f_{23},e_2,e_3$ which has the same union but they have two private parts: the vertices of color $4$ in $e_2$ and the vertex of color $1$ in $e_3$.
This reduces Case 3 to Case 2.

\noindent {\bf Case 4.} None of the edges $e_1,e_2,e_3$ have private parts.

\noindent Vertices uncovered by $X$ are colored with $1$. Here $f_{12},f_{13},f_{14},f_{15},f_{23},f_{24},f_{25}$ are all excluded since there is some $e_i$ intersecting them in at most one vertex. The other three edges $f_{34},f_{35},f_{45}$ are excluded since they are proper subsets of some $e_i$.

In all cases we found a $3$-strong coloring with at most five colors. \qed

\end{document}